\documentclass[12pt]{article}
\newcommand{\qed}{\ifmmode$\Box$\else{\unskip\nobreak\hfil
\penalty50\hskip1em\null\nobreak\hfil$\Box$
\parfillskip=0pt\finalhyphendemerits=0\endgraf}\fi}
\usepackage{amsmath}
\usepackage{latexsym}
\usepackage{amsfonts}
\begin{document}
\renewcommand{\theequation}{\thesection.\arabic{equation}}
\newtheorem{Pa}{Paper}[section]
\newtheorem{Tm}[Pa]{{\bf Theorem}}
\newtheorem{La}[Pa]{{\bf Lemma}}
\newtheorem{Cy}[Pa]{{\bf Corollary}}
\newtheorem{Rk}[Pa]{{\bf Remark}}
\newtheorem{Pn}[Pa]{{\bf Proposition}}
\newtheorem{Ee}[Pa]{{\bf Example}}
\newtheorem{Dn}[Pa]{{\bf Definition}}
\newtheorem{St}[Pa]{{\bf Statement}}
\newtheorem{I}[Pa]{{\bf}}
\oddsidemargin 0in

{\hfill IFT UwB  03/2002}

\bigskip

\begin{center}
{\bf \Large On the structure of $C^*-$algebras generated by the components of
polar decompositions\footnote{This work is supported in part by KBN grant 2 PO3 A 012 19.}}\\
\bigskip
{\large A. Lebedev,   A. Odzijewicz  }\\ \bigskip
 Belarus State University,
Department of Mathematics and Mechanics, Skariny av. 4, Minsk,
220050, Belarus;\\

Institute of Theoretical Physics,  University  in Bialystok,
 ul. Lipowa 41, PL-15-424 Bialystok, Poland
\end{center}
\bigskip
\begin{abstract}
In the present paper we study the structure of $C^*-$algebras generated by the components of
the polar decompositions of  operators in Hilbert space satisfying certain commutation relations.
\end{abstract}



\tableofcontents
\section{Introduction}
\setcounter{equation}{0}
Let $a$ be a linear
bounded operator acting in a Hilbert space $H$. Suppose that the operator $a$
and its adjoint $a^*$ satisfy the following  relation
\begin{equation}
\label{eI.1}
aa^* \in \{1, a^*a \},
\end{equation}
where we denote by $ \{ \alpha_\lambda , \  \lambda\in\Lambda \}$
 the $C^*-$algebra generated by the family of operators
$\alpha_\lambda , \  \lambda\in\Lambda $.
The problems we shall deal with in this paper concern the  description of
 the structure of the $C^*-$algebra $\{ 1, a \}$.\\
\begin{I}\em
$C^*-$algebras defined by  condition (\ref{eI.1}) are the natural
generalizations of algebras satisfying the relation
\begin{equation}
\label{e.gamma} aa^* = \gamma (a^* a)
\end{equation}
where $\gamma$ is a continuous function on the spectrum of the
operator $a^*a$. they appear in different problems of quantum
physics (see, for example, \cite{Kar,OdHorTer}).\\

\noindent In the simplest case when $\gamma $ is a linear map
\begin{equation}
\label{e.q} aa^* = q\, a^* a + h
\end{equation}
they include: \\

(i) the  commutative algebras generated by normal operators:
$$a^*a = aa^*$$  ($q=1, \ \ h=0$)\ \ the theory of which  forms a
cornerstone of the spectral theory of operators in Hilbert
spaces.\\

(ii) Toeplitz  algebra: $$aa^* =1 $$ ($q=0, \ \ h=1$)\\

(iii) Heisenberg algebra:
\begin{equation}
\label{e.heisen}  aa^* - a^*a =h
\end{equation}
 ($q=1, \ \ h>0$).\\ This is in fact {\em not} a $C^*-$algebra
 situation: the operators satisfying (\ref{e.heisen}) are
{\em unbounded} and play the  principal role  (associated with the
position and momentum operators)
 in the quantum mechanics theory.\\

\noindent In the general case \ \ ($q\neq 0,1 \ \ h\neq 0$)\ \ the quantum algebra generated by
relation (\ref{e.q}) could be interpreted as a $q-$deformation of the Heisenberg algebra. The
models of quantum mechanics based on the $q-$deformation of the Heisenberg algebra has been
studied, for example, in \cite{MaxOdz}.\\

\noindent In the physics of many bosons system that appear in a
natural way in quantum optics and nuclear physics one arrives at
the Shrodinger operator $H$ given by
\begin{equation}
\label{e.H} H = a + a^* + D (a^* a)
\end{equation}
where the operator $a$ satisfies (\ref{eI.1}) (see
\cite{OdHorTer}).\\ For the corresponding quantum system the
algebra arising from (\ref{eI.1}) plays the role of 'symmetry
algebra'. It enables us to establish the interrelation between the
spectral problem for the operators considered and the theory of
orthogonal polynomials and to solve the problem in many concrete
cases.\\

\noindent Let us mention also  that the algebra given by
(\ref{eI.1})  finds applications in the theory of basic
hypergeometric functions and integrable systems (see
\cite{Odz1}).\\

\noindent Along with the algebra $\{1,a  \}$ it is reasonable to
consider the following objects and arising problems.\\ Let
\begin{equation}
\label{eI.2}
a =U | a |
\end{equation}
be the standard polar decomposition of $a$. Here\ \  $ | a | = \sqrt{a^*a}$\ \ and \  $U$ \  is
a partial isometry
defined by
\begin{equation}
\label{eI.3}
U(|a|\xi ) = a\xi , \ \ \ \xi \in H .
\end{equation}
Then (\ref{eI.1}) means that
\begin{equation}
\label{eI.4}
U |a|^2 U^* \in \{ 1,  |a|^2 \} =  \{ 1,  |a| \}.
\end{equation}
Note also that since $U^*U$ is the orthogonal projection onto ${\rm Im}\, |a|$ it follows
that $U^*U |a| =|a|$ and therefore
\begin{equation}
\label{eI.5}
U |a|^2 U^* = U |a| U^*U |a| U^* = (U |a| U^*)^2 .
\end{equation}
Thus we conclude that (\ref{eI.1}) is equivalent to
\begin{equation}
\label{I.a}
 U |a| U^* \in \{1, |a| \} .
\end{equation}
In addition (\ref{eI.5}) means  that the mapping
 $$
 U(\cdot )U^* \, :  \{  |a| \} \, \to \, \{ 1, |a| \}
 $$
 is  a {\em morphism}.\\
 This leads to the tight interrelation of the algebras induced by (\ref{eI.1})   with the crossed products of $C^*-$algebras by
 semigroups of endomorphisms. This interrelation is clarified in
 section 4 (see in particular Theorem \ref{3a.5} and Remark
 \ref{cross}).\\

 \noindent  In connection with the matter considered in the paper it is reasonable to mention
 that a great number of various irreducible representations of the relation
 (\ref{eI.1}) along with the physical applications of the $C^*-$algebras
 arising are presented in the book \cite{OstSam}.\\
\end{I}

\begin{I}
\label{content} \em The main aim of  the article is to reveal the
structure  of the $C^*-$algebras
  $\{ 1, |a|, U \}$  and $\{1, a  \} \subset \{ 1, |a|, U \}$ \ \ where $a$ satisfies relation
 (\ref{eI.1}) and $U$ is the partial isometry taken from the polar decomposition (\ref{eI.2}) . \\

\noindent  The paper is organized as follows. \\ We start with an
auxiliary section \ref{1} where we present a number of (mostly
known)  facts on partial isometries and the associated mappings of
the algebra of linear operators in a Hilbert space that we shall
use in the subsequent sections.\\

\noindent Sections \ref{2} and \ref{3} present the main framework
(the general operator algebraic picture that one can attend 'from
on high') of the structures and objects that should be involved in
the investigation of the algebras in question.\\ As a result we
obtain the necessary 'coefficient' algebra (given in Theorem
\ref{A.inf}) and the crucial property \ref{*} that enables us to
reveal the structure of the arising algebras up to $^*-$algebraic
isomorphism (Theorem \ref{3a.5}).\\

\noindent Finally in sections \ref{4} and \ref{5} we apply the
results obtained to the initial algebras \ \ $\{ 1, |a|, U \}$ and
$\{1, a  \}$\ \  under investigation and give the description of
their structure.
\end{I}

\section{Preliminaries. Partial isometries and endomorphisms}
\label{1}
\setcounter{equation}{0}

This is a starting  auxiliary section. Here we list some facts
about partial isometries and the associated mappings  of the
algebra of linear operators in a Hilbert space that we shall use
in the subsequent sections. Most of the facts presented here are
known and we give the proofs for the sake of completeness.\\

 \noindent Recall that a
linear bounded operator $U$ in a Hilbert space $H$ is called a
partial isometry if there exists a closed subspace $H_1 \subset H$
such that $$ \Vert U\xi \Vert = \Vert \xi  \Vert , \ \ \ \xi \in
H_1 $$ and $$ U\xi = 0, \ \ \ \ \xi \in H\ominus H_1. $$ The space
$H_1$ is called the {\em initial} space of $U$ and $U(H_1)$ is
called the {\em final} space of $U$.\\

\noindent Hereafter we list the  well known characteristic
properties of a partial isometry  (see for example \cite{Halm},
problem 98):\\
\begin{Tm}
\label{2.1}
The following statements are equivalent:\\

\noindent 1) $U$ is a partial isometry,\\
2) $U^*$ is a partial isometry,\\
3) $U^*U$ is a projection (onto the initial space of $U$),\\
4) $UU^*$ is a projection (onto the final space of $U$),\\
5) $UU^*U =U$ and \ $U^*UU^* =U^*$.\\
\end{Tm}

\noindent One can associate with any linear bounded operator  $V$
acting in a Hilbert space $H$ the following  mapping $\delta_V :
L(H)\to L(H) $ of the algebra $L(H)$ of all linear bounded
operators in $H$:
\begin{equation}
\label{e2.1} \delta_V (b) = VbV^*, \ \ \  b\in L(H)
\end{equation}
It is clear that $\delta_V$ is a linear positive mapping and
$\delta _V (b^*) = [\delta_V (b)]^*$ \ \  for any $b\in L(H)$.\\

\noindent The next proposition shows that $\delta_V$ generates an
{\em endomorphism} of some subalgebra $A\subset L(H)$ only under
rather restrictive  assumptions on $V$ (namely $V$ should be a
partial isometry  possessing special additional properties).
\begin{Pn}
\label{endom} For a given operator $V\in L(H)$ the following
statements are equivalent:\\

\noindent (i) There exists a subalgebra $A\subset L(H)$ containing
the identity such that the mapping $$ V (\cdot ) V^* : L(H) \to
L(H) $$ generates an endomorphism of $A$.\\

\noindent (ii) All the operators $V^k, \ k=1,2, ...$ are partial
isometries.\\

\noindent (iii) The family $V^{*k}V^k, \ k=1,2,...$ is a
commutative decreasing family of projections.
\end{Pn}
{\bf Proof:} In view of statement 3) of Theorem \ref{2.1} (iii)
implies (ii).\\ On the other hand if (ii) is true then
$V^kV^{*k}V^k =V^k, \ \ k=1,2, ...$ by statement 5) of Theorem
\ref{2.1}. Therefore for any $k\ge l$ we have $$
V^{*k}V^kV^{*l}V^l = V^{*k}V^{k-l}V^lV^{*l}V^l = V^{*k}V^{k-l}V^l
= V^{*k}V^k$$ and in the same way one can verify that
$$V^{*l}V^lV^{*k}V^k =  V^{*k}V^k .$$ Thus (ii) implies (iii).\\
  If (i) is true then all
the mappings\ \  $V^k (\cdot )V^{k*}, \ k=1,2,...$\ \ are also
endomorphisms of the algebra $A$ (mentioned in (i)) and therefore
(since $1\in A$) $$ V^kV^{k*}=V^k 1V^{k*} = V^k (1\cdot
1)V^{k*}=(V^k1V^{k*})(V^k1V^{k*})= (V^kV^{k*})^2 \ , \ \ \ k=1,2,
... $$ Thus (ii) is true.\\ Finally if (iii) is true then\ \
$V(\cdot )V^*$\ \ generates an endomorphism of the algebra \ \ $A
= \{ 1, VV^*, V^2V^{2*}, ... , V^k V^{k*}, ... \}$.
\mbox{}\qed\mbox{}\\

\noindent In view of this result it is reasonable to consider the
mapping $\delta_V$ in the situation when $V$ is a partial isometry
in more detail.\\ In the next lemma we list some  simple
properties of the mapping $\delta_V$ in this case.
\begin{La}
\label{2.3} Let $V$ be a partial isometry  then\\

\noindent 1) $\delta_V \delta_{V^*}\delta_V = \delta_V$\ \ and \ \
$\delta_{V^*}\delta_V \delta_{V^*} =\delta_{V^*}$\\

\noindent 2) For any $b\in L(H)$ we have
 $$\delta_V (b) = VV^* \delta_V (b) = \delta_V (b)VV^* \ \ { and}\ \  \
\delta_{V^*}(b) = V^*V \delta_{V^*}(b) = \delta_{V^*}(b)V^*V .$$

\noindent 3) ${\rm Ker}\, \delta_V  = \{ b\in L(H): V^*V b V^*V =0
\}$.\\ In particular if \ \ $V^*V =1$ (that is $V$ is an isometry)
then ${\rm Ker}\, \delta_V = 0$\\

 \noindent 4) ${\rm Ker}\, \delta_{V^*} = \{ b\in L(H): VV^* bVV^* =0 \}$.\\
In particular if \ \ $VV^* =1$ (that is $V^*$ is an isometry) then
${\rm Ker}\,\delta_{V^*} = 0$\\
\end{La}
{\bf Proof:} The first and the second statements follow from the
equalities \ \ $VV^*V =V$ \ \  and \ \ $V^* VV^* =V^*$.\\ Let us
verify the third statement. \\ If\ \  $ V^*V b V^*V =0$\ \  then
$$
 \delta_V (b) = VbV^* = VV^*V b V^*VV^* =  V(V^*V b V^*V)V^* =0
$$ that is $b\in {\rm Ker}\, \delta_V $.\\ On the other hand if\ \
$b\in {\rm Ker}\, \delta_V $ \ \ then $$ 0 = VbV^* = V^*(VbV^*)V =
V^*VbV^*V $$ Thus 3) is proved.\\

\noindent Statement 4) can be proved in the same way.
\mbox{}\qed\mbox{}\\

\noindent In fact if we do not look for endomorphisms of some
algebra $A\subset L(H)$ (as in Proposition \ref{endom}) but only
for {\em morphisms} \ $\delta_V :A\to \delta_V (A)$ nevertheless
(as the next lemma shows) this situation is closely related to the
property of $V$ being a partial isometry.
\begin{La}
\label{V} Let  $A$ be a $C^*-$subalgebra of $L(H)$.\\

 {\bf I}.\ \ If\ \ $V$ is a partial isometry and \ \  $V^* V \in
A^\prime$\ \ (where $A^\prime$ is the commutant of $A$) then\\

\noindent (i) $\delta_V :A\to \delta_V (A)$ is a morphism,\\

\noindent (ii) $V\alpha = \delta_V (\alpha ) V$\ \ and \ \ $\alpha
V^* =V^* \delta_V (\alpha ) $\ \  for every $\alpha \in A$.\\

{\bf II}.\ \ If $A$  contains the identity and \ \  $VV^* , V^* V
\in A^\prime$\ \ then the following statements are equivalent:\\

\noindent (i)\ \  $V^*V$ is a projection (that is $V$ is a partial
isometry),\\

\noindent (ii)\ \  $VV^*$ is a projection,

\noindent \\ (iii)\ \ $V(\cdot )V^* : A \to VAV^*$ is a
morphism,\\

\noindent (iv)\ \  $V^* (\cdot )V : A \to V^* A V$ is a morphism.
\end{La}
{\bf Proof:}\ \ {\bf I}. The first statement here follows from the
observation that for any \ \ $\alpha , \beta \in A$ \ \ we have $$
V(\alpha \cdot \beta )V^* = VV^*V \alpha \cdot \beta V^* =V \alpha
V^*V \beta V^*= (V  \alpha V^*)(V  \beta V^*) $$ and the second
follows from the relations $$ V\alpha = VV^*V \alpha = V \alpha
V^*V = (V\alpha V^*) V $$ along with the passage to the adjoint
operators.\\

\noindent {\bf II}. Here (i) and (ii) are equivalent in view of
Theorem \ref{2.1}.\\ If (i) is true then $V$ is a partial isometry
and (iii) follows from statement (i) in {\bf I}.\, Similarly (ii)
implies (iv). On the other hand if (iii) is true then $$ VV^*=V
1V^* = V (1\cdot 1)V^*=(V1V^*)(V1V^*)= (VV^*)^2 $$ thus (ii) is
true.\\ Similarly (iv) implies (i).
 \mbox{}\qed\mbox{}

\noindent The next result pushes the matter a bit further and
clarifies that if in addition to what is presumed in part  I of
Lemma \ref{V} \ \ $\delta _V$ generates an {\em endomorphism} of
$A$ then the interrelation between $A$ and $V$ is much closer.
\begin{La}
\label{endom1} Let $V$ be a partial isometry and $A$ be a
$C^*-$subalgebra of $L(H)$.
 If\ \  $V^* V \in A^\prime$\ \  and\ \  $VAV^* \subset A$\ \  then\\

\noindent $\delta_V : A \to A$ is an endomorphism and\ \
$V^{*k}V^k \in A^\prime \ , \ \ k=1,2, ...$
\end{La}
{\bf Proof:} By part I of Lemma \ref{V} in the situation under
consideration we have that\ \ $\delta_V : A \to A$\ \ is an
endomorphism.\\ Applying statement (ii) of part I of Lemma
\ref{V}\ \ $k$ times we have for any $\alpha \in A$ $$ V^{*k}V^k
\alpha = V^{*k}(V^k \alpha V^{*k}) V^k $$ and $$
 \alpha V^{*k}V^k =  V^{*k}(V^k \alpha V^{*k}) V^k
$$ Thus \ \ $V^{*k}V^k \in A^\prime \ , \ \ k=1,2, ...$.
\mbox{}\qed\mbox{}\\
 In connection with the foregoing  lemma it is
reasonable to mention the following  properties of partial
isometries
\begin{La}
\label{""} Let $V$ be a partial isometry such that
\begin{equation}
\label{e.000} [V^*V, V^kV^{*k}]=0 \ , \ \ \ k=1,2, ...
\end{equation}
where $[\alpha ,\beta ]$ is the commutator of the sets $\alpha $ and\ \ $\beta$. Then\\

\noindent 1)\ \
$
V^{*l}V^l \in \{ VV^*, V^2V^{*2}, . . . , V^kV^{*k}, . . .
\}^{\prime}\ , \ \ l=1,2,... $\\

\noindent 2) \ \ For any $1\le k\le l $ \ \  $V^* V^k V^{*l} =
V^{k-1}V^{*l}$\\

\noindent 3)\ \ All the operators $V^k$ are partial isometries and
the  family $V^kV^{k*}, \ k=1,2,...$ is a commutative decreasing
family of projections.
\end{La}
{\bf Proof:} 1) This follows from Lemma \ref{endom1} since   the
mapping $V(\cdot )V^*$
 is an endomorphism of the algebra \ \
$ \{ VV^*, V^2V^{*2}, . . . , V^kV^{*k}, . . . \}$.\\

\noindent 2). Let\ \ $1\le k\le l $\ \ then applying
(\ref{e.000}) we obtain $$
\begin{array}{*{20}l}
V^* V^k V^{*l} = (V^*V) (V^{k-1}V^{*k-1})V^{*l-k+1}=
(V^{k-1}V^{*k-1})(V^*V)V^{*l-k+1} =\\ \ \\ (V^{k-1}V^{*k-1})
(V^*VV^*) V^{*l-k} = V^{k-1}V^{*l}
\end{array}
$$

\noindent 3).  Applying 2)\ \  \mbox{$k$\  times} we have $$
\begin{array}{*{20}l}
(V^kV^{*k})^2 = V^kV^{*k}V^kV^{*k}= V^kV^{*k-1}(V^*V^kV^{*k} ) =
V^kV^{*k-1}V^{k-1}V^{*k} = \\ \ \\ V^k V^{*k-2}(V^*
V^{k-1}V^{*k-1})V^* = V^k V^{*k-2}V^{k-2}V^{*k} = ... = V^kV^{*k}
\end{array}
$$ this means that $V^k \ , \ \ k=1,2, . . . $ are isometries and
therefore  3) is true. \mbox{}\qed\mbox{}\\



\section{Extensions of commutative $C^*-$algebras by endomorphisms }
\label{2}
 \setcounter{equation}{0}

 Throughout this section we fix a partial isometry
$U$ and a commutative $C^*-$algebra
$A_0\subset L(H)$ containing the identity and
 proceed to the description of the structures and relations associated with
the mappings\ \ $\delta_U$\ \ and\ \ $\delta_{U^*}$
(\ref{e2.1}).\\ To shorten the notation we denote\ \ $\delta_U$\ \
simply by\ \ $\delta$\ \ and\ \ $\delta_{U^*}$\ \ by\ \
$\delta_*$\ \ thus
\begin{equation}
\label{e.delta}
 \delta (b)= UbU^*, \ \ \ b\in L(H)
\end{equation}
and
\begin{equation}
\label{e.delta*} \delta_* (b)= U^*bU, \ \ \ b\in L(H)
\end{equation}\\

\noindent Our goal in this section is the following:\\

starting from the algebra $A_0$ and the mappings\ \ $\delta$\ \
and \ \ $ \delta_*$\ \ we wish to construct  a {\em commutative}
extension\ \ $A\supset A_0$\ \ such that both the mappings \ \
$\delta$\ \ and \ \ $ \delta_*$\ \ are endomorphisms of $A$.\\

\noindent We shall construct the algebra $A$ in two steps. At the
first step we construct the extensions $A_\infty$ \ \ ($_\infty
A$) \ \ (see (\ref{e.A1}) and  (\ref{e.A2})) such that \ \
$\delta$ \ \ (or $\delta_*$ \ \ respectively) generate their
endomorphisms. This is done in Theorems \ref{A1} and \ref{A2}.
Then as the second step we extend the algebra\ \ $A_\infty$ \ \
($_\infty A$) \ \ further to obtain the algebra \ \ $A =\, _\infty
(A_\infty)\, = (\, _\infty A)_\infty$ \ \ (see (\ref{e.inf}) and
(\ref{e.Ainfty2})) that is stable as with respect to \ \ $\delta$\
\ so also with respect to \ \ $\delta_*$. This is done in Theorems
\ref{A0.inf} and \ref{A.inf} which also describe the subtle
internal hierarchical structure of the algebras arising in the
process of extension.\\

\noindent It is reasonable to note that since we are looking for a
{\em commutative} extension $A$ of $A_0$ and in view of Lemmas
\ref{V} and \ref{endom1} among the natural assumptions on \ \
$\delta$ \ \ (or $\delta_*$ \ \ respectively) we have \ \
$\delta^k (1) \in A_0^\prime$ \ \ ( $\delta_*^k (1) \in
A_0^\prime$).\\

\noindent Now we start to produce the desired extensions of
$A_0$.\\
\begin{Tm}
\label{commutant} Let \ \ $\delta_*^k (1), \ k=\overline{0,n}$ \ \
be the projections and\ \
 $\delta_*^k (1) \in A_0^\prime$ \ \ and
  $\delta^k (A_0) \subset A_0^\prime , \ k=\overline{0,n}$\ \
(where $\delta_*^0 (x) = \delta^0 (x) = x, \ \ x\in L(H)$, in
particular $\delta_*^0 (1) = \delta^0 (1) = 1$).
 Then
\begin{equation}
\label{e.com1} [\delta_*^k (A_0), \delta_*^l (A_0)] =0 , \ \ \
0\le k,l \le n
\end{equation}
and
\begin{equation}
\label{e.com2} [\delta^k (A_0), \delta^l (A_0)] =0 , \ \ \ 0\le
k,l \le n
\end{equation}
\end{Tm}
{\bf Proof:}
Take some $ k\ge l$ then for any $\alpha , \beta \in A_0$ we have
$$
\begin{array}{*{20}l}
 {\delta_*^k (\alpha )\cdot \delta_*^l(\beta ) = U^{*k}\alpha U^k U^{*l}\beta U^l =
U^{*k}\alpha U^{k-l}(U^l U^{*l})\beta U^l =} \\
{\ \ }\\
{U^{*k}\alpha U^{k-l}\beta (U^l U^{*l})U^l = U^{*k}\alpha U^{k-l}\beta U^l =
 U^{*k}\alpha U^{k-l}(U^{*k-l}U^{k-l})\beta U^l =}\\
{\ \ }\\ {U^{*k}\alpha U^{k-l}\beta U^{*k-l}U^{k-l}U^l =
U^{*k}\alpha (U^{k-l}\beta U^{*k-l})U^k =}\\ {\ \ }\\
{U^{*k}\alpha  \delta^{k-l} (\beta )    U^k = \delta_*^k (\alpha
\cdot \delta^{k-l} (\beta ))}
\end{array}
$$
On the other hand
$$
\begin{array}{*{20}l}
{ \delta_*^l(\beta )\cdot \delta_*^k (\alpha ) = U^{*l}\beta U^l
U^{*k}\alpha U^k = U^{*l}\beta (U^lU^{*l}) U^{*k-l}\alpha U^k =
}\\ {\ \ }\\ { U^{*l}U^lU^{*l}\beta U^{*k-l}\alpha U^k =
U^{*l}\beta  U^{*k-l}\alpha U^k = U^{*k} (\delta^{k-l} (\beta
)\cdot \alpha  )U^k =}\\ {\ }\\ {\delta_*^k (\delta^{k-l} (\beta
)\cdot \alpha  )}
\end{array}
$$ Comparing these two relations and bearing in mind that by
assumption of the theorem \ \ $\delta^{k-l} (\beta )\in
A_0^\prime$ \ \ we conclude that $$ \delta_*^k (\alpha \cdot
\delta^{k-l} (\beta ))=\delta_*^k (\delta^{k-l} (\beta )\cdot
\alpha  ) $$ and therefore (\ref{e.com1}) is proved.\\ To prove
(\ref{e.com2}) we take \ \ $k,\ l,\ \alpha , \ \beta $ \ \ as
above and observe that $$
\begin{array}{*{20}l}
{\delta^k (\alpha )\cdot \delta^l(\beta )= U^k \alpha U^{*k}U^l
\beta U^{*l}= U^k \alpha U^{*k-l}(U^{*l}U^l) \beta U^{*l}=}\\ {\
}\\ { U^k \alpha U^{*k-l} \beta U^{*l}U^lU^{*l}= U^k \alpha
U^{*k-l} \beta U^{*l}=
 U^l (U^{k-l} \alpha U^{*k-l} )\beta U^{*l}= }\\
{\ }\\ {U^l ( \delta^{k-l} (\alpha )\cdot \beta )U^{*l}= \delta^l
( \delta^{k-l} (\alpha )\cdot \beta )}
\end{array}
$$
On the other hand
$$
\begin{array}{*{20}l}
{\delta^l(\beta )\cdot \delta^k (\alpha ) = U^l \beta U^{*l}U^k
\alpha U^{*k} = U^l \beta (U^{*l}U^l)U^{k-l} \alpha U^{*k} = }\\
{\ }\\ {U^l U^{*l}U^l\beta U^{k-l} \alpha U^{*k} =U^l \beta
U^{k-l} \alpha U^{*k} = U^l (\beta  \cdot \delta^{k-l} (\alpha )
)U^{*l} =}\\ {\ }\\ {\delta^l (\beta  \cdot \delta^{k-l} (\alpha )
) }
\end{array}
$$
which (as in the previous situation) implies (\ref{e.com2}).
\mbox{}\qed\mbox{}
\begin{I}
\label{defAn}
\em
For every $n= 0,1, ...$ set
\begin{equation}
\label{e.An1} A_n  = \{ A_0, ..., \delta^n (A_0) \} .
\end{equation}
and
\begin{equation}
\label{e.An} _nA  = \{ A_0, ..., \delta_*^n (A_0) \} .
\end{equation}
\end {I}
\begin{Tm}
\label{An} Let\ \ $\delta_*^k (1)\ , \ \ k=\overline{1,n}$\ \ be
the projections and $$ \delta (A_0)\subset A_0 \ \ \ and \ \ \
\delta_*(1)\in {A_0}^\prime $$
 then\ \  $_nA $\ \ is a commutative algebra and \\

\noindent (i) $\delta_*^k (1) \in {A_0}^\prime \ , \ \
k=1,2,...$\\

\noindent (ii) For any\ \  $0 \le l \le k \le n$ \ \  we have \ \
$$ \delta_*^k (A_0)\cdot \delta_*^l (A_0)  \subset \delta_*^k
(A_0) $$ in particular\ \  $\delta_*^n (A_0) $\ \ is an ideal in\
\ $_nA $,\\

\noindent (iii)\ \ $_nA $ \ \ is the set of  operators of the form
$$ \alpha_0 + \delta_*(\alpha_1) + ... + \delta_*^n(\alpha_n), \ \
\alpha_i \in A_0. $$
 (iv)
$\delta  :\,   _nA \,  \to \ \ _{n-1}A $\ \  is a morphism.
\end{Tm}
{\bf Proof:} Since $\delta (A_0)\subset A_0 $ \ \ and \ \ $
\delta_*(1)\in {A_0}^\prime$\ \ Lemma \ref{endom1} implies (i).\\
Therefore\ \ (as $\delta^k (A_0)\subset A_0 \subset {A_0}^\prime \
, \ \ k=1,2,...$)\ \ all
  the conditions of Theorem  \ref{commutant}
are satisfied   it follows that\ \  $_nA $\ \ is commutative.\\ To
verify (ii) observe that in the proof of Theorem \ref{commutant}
it was established  that for any \ $k\ge l$\  and \ $\alpha ,\beta
\in A_0$ \ we have $$ \delta_*^k (\alpha )\cdot \delta_*^l(\beta )
=
 \delta_*^k (\alpha \cdot \delta^{k-l} (\beta ))
$$ Since by the assumption we have\ \ $ \delta^{k-l} (\beta ) \in
A_0$ \ \ this equality implies (ii).\\ Clearly (iii) follows from
(ii).\\ Finally  by the assumption we have\ \ $UU^* = \delta(1)
\in A_0$.\  Therefore we conclude (using (ii)) that for any
$\alpha \in A_0$ and any $1\le k\le n$ the following relations
hold $$
 \delta(\delta_*^k (\alpha )) = UU^{*k} \alpha U^k U^* = UU^* \delta_*^{k-1} (\alpha ) UU^*=
\delta_*^{k-1} (\alpha ) UU^* \in \delta_*^{k-1}(A_0). $$ Thus\ \ $\delta \, :\, _nA
\hookrightarrow \, _{n-1}A $.\\ Note also that $$ U^*U = \delta_*(1) \in\, _nA \subset (_nA
)^\prime $$ and therefore $\delta \, :\, _nA \hookrightarrow \, _{n-1}A $ is a morphism by Lemma
\ref{V}. So (iv) is proved. \mbox{}\qed\mbox{}\\ The next two theorems describe the natural
extensions of the algebra $A_0$ that are stable under the morphisms $\delta$ and $\delta_*$
respectively.
\begin{Tm}
\label{A1} Let \ \ $\delta_*^k (1), \ k= 0, 1, ...$ \ \ be the
projections and $$ \delta_*^k (1) \in A_0^\prime , \ \ \ \delta^k
(A_0) \subset A_0^\prime ,\ \ \ \delta_* (1) \in [\delta^k
(A_0)]^\prime , \ \  k= 0, 1, ... $$ then
\begin{equation}
\label{e.A1}  A_\infty = \{ A_0, \delta (A_0), . . . , \delta^k
(A_0), . . .  \}
\end{equation}
is the minimal commutative $C^*-$algebra containing $A_0$ and such
that $$ \delta : A_\infty  \to A_\infty  $$ is an endomorphism.\\
In addition we have
\begin{equation}
\label{e.delta1} U\alpha = \delta(\alpha )U\ \ \ \mbox{for any}\ \
\alpha \in A_\infty
\end{equation}
and
\begin{equation}
\label{e.com} \delta_*^k (1) \in (A_\infty  )^\prime \ , \ \ \
k=1,2,...
\end{equation}
\end{Tm}
{\bf Proof:} The commutativity of $A_\infty $ follows from Theorem
\ref{commutant}.\\ The minimality of  $A_\infty $ is clear from
the construction.\\ So let us  verify  that \
$
\delta $\ is an endomorphism of $A_\infty $.\\ Take any \ \
$k,l\ge 0$\ \  and \ \ $\alpha , \beta \in A_0$.  We have $$
\begin{array}{*{20}l}
{\delta   (\delta^k (\alpha )\cdot  \delta^l (\beta )) = U
\delta^k (\alpha ) \delta^l (\beta ) U^*= U (U^*U) \delta^k
(\alpha ) \delta^l (\beta ) U^*= }\\ {\ }\\ {U \delta^k (\alpha )
U^*U \delta^l (\beta ) U^*= \delta(\delta^k (\alpha )) \cdot
\delta (\delta^l (\beta ))}
\end{array}
$$ which implies that  \ $\delta$ \ is an endomorphism of
$A_\infty $.\\ In addition the assumption $$ \delta_* (1) \in
[\delta^k (A_0)]^\prime , \ \  k= 0, 1, ... $$ implies\ \
$\delta_* (1) \in  (A_\infty  )^\prime$\ \ and therefore
(\ref{e.delta1}) follows from Lemma \ref{V}.\\ Finally since\ \
$\delta_* (1) \in (A_\infty  )^\prime$\ \ and   \ $\delta$ \ is an
endomorphism of $A_\infty $\ \ (\ref{e.com}) follows from Lemma
\ref{endom1}. \mbox{}\qed\mbox{}
\begin{Tm}
\label{A2} Let\ \  $\delta_*^k (1), \ k= 0, 1, ...$\ \  be the
projections and $$ \delta_*^k (1) \in A_0^\prime , \ \ \ \delta^k
(A_0) \subset A_0^\prime ,\ \ \ \delta (1) \in [\delta_*^k
(A_0)]^\prime , \ \  k= 0, 1, ... $$ then
\begin{equation}
\label{e.A2} _\infty A  = \{ A_0, \delta_* (A_0), . . . ,
\delta_*^k (A_0), . . .  \}
\end{equation}
is  the minimal commutative $C^*-$algebra containing $A_0$ and
such that $$ \delta_* \, :\,  _\infty A \,   \to \,  _\infty A  $$
is an endomorphism.\\ In addition  we have
\begin{equation}
\label{e.antidelta} U^*\alpha = \delta_*(\alpha )U^*\ \ \
\mbox{for any}\  \ \ \alpha \in \, _\infty A
\end{equation}
and
\begin{equation}
\label{e.antidelta1} \delta^k (1) \in (_\infty A  )^\prime \ , \ \
\ k=1,2,...
\end{equation}
\end{Tm}
{\bf Proof:} The commutativity of\ \ $_\infty A $\ \ follows from
Theorem \ref{commutant}.\\ The minimality of\ \  $_\infty A $\ \
is clear from the construction.\\ The property that $\delta_*$ is
an endomorphism of\ \ $_\infty A $\ \ follows from the observation
that for any\ \  $k,l \ge 0$\ \ and \ \ $\alpha  , \beta \in A_0$\
\ we have $$
\begin{array}{*{20}l}
{\delta_* (\delta_*^k (\alpha ) \cdot \delta_*^l (\beta )) = U^*
\delta_*^k (\alpha ) \delta_*^l (\beta ) U= U^* (UU^*) \delta_*^k
(\alpha ) \delta_*^l (\beta ) U= }\\ {\ }\\ {U^* \delta_*^k
(\alpha )UU^* \delta_*^l (\beta ) U=
 \delta_* (\delta_*^k (\alpha )) \cdot \delta_* (\delta_*^l (\beta )) }
\end{array}
$$ In addition the assumption $$ \delta (1) \in [\delta_*^k
(A_0)]^\prime , \ \  k= 0, 1, ... $$ implies\ \ $\delta (1) \in
(_\infty A    )^\prime$\ \ and thus (\ref{e.antidelta}) follows
from Lemma \ref{V}.\\ Finally (\ref{e.antidelta1}) can be
established just in the same way as it was done for (\ref{e.com})
in the proof of the previous theorem. \mbox{}\qed\mbox{}
\begin{Rk}
\label{asymm} \em One can notice certain asymmetry  (between\ \
$\delta$\ \ and\ \ $\delta_*$) as in the statement  of Theorem
\ref{commutant} so also in the statements of Theorems \ref{A1} and
\ref{A2}. This means that we can also formulate and prove the
analogous statements by exchanging the conditions on\ \ $\delta$\
\ and\ \ $\delta_*$.
\end{Rk}

\noindent Now we present the description of the extensions of the
algebra\ \ $A_0$\ \ that are stable as with respect to\ \
$\delta$\ \ so also with respect to\ \ $\delta_*$.
\begin{Tm}
\label{A0.inf} Let\ \ $\delta (A_0) \subset A_0 $.\\

 {\bf I}.\ \
The following statements are equivalent:\\

\noindent (i) \ \ There exists a commutative $C^*-$algebra\ \
${\cal A } \supset A_0$\ \ such that both the mappings \ \
$\delta$\ \ and \ \ $\delta_*$ are endomorphisms of \ \ ${\cal
A}$.\\

\noindent (ii)\ \ All the operators \ \ $\delta_*^k (1), \ k= 0,
1, ...$\ \ are projections and\ \ $\delta_* (1) \in  A_0^\prime
$.\\

{\bf II}. \ \ If condition (ii) of part {\bf I} is satisfied then
 the algebra\ \ $_\infty A  $\ \ (\ref{e.A2})
is the minimal commutative $C^*-$algebra such that\ \ $\delta$\ \
and\ \ $\delta_*$\ \ are  endomorphisms of\ \ $_\infty A   $.
\\
 Moreover we have:\\

\noindent (i) For any\ \  $0 \le l \le k $ $$ \delta_*^k
(A_0)\cdot \delta_*^l (A_0)  \subset \delta_*^k (A_0) . $$

\noindent (ii) For every $n$ the subalgebra\ \  $_nA \subset \,
_\infty A  $\ \ defined by (\ref{e.An}) is the set of the
operators of the form $$
 \alpha_0 + \delta_*(\alpha_1) + ... +
\delta_*^n(\alpha_n), \ \ \alpha_i \in A_0. $$

\noindent (iii) $\delta \, :\, _nA  \to \,  _{n-1}A, \ \ n=1,2,
... $ \ \ \ \ \
      $\delta_* \, :\, _nA \,  \to \, _{n+1}A, \ \ n=0,1, ... $,\\

\noindent (iv)
$
U\alpha = \delta (\alpha )U \ \  and \ \ \ U^* \alpha
=\delta_*(\alpha )U^* $\ \ for any \ \ $\alpha \in \, _\infty A $.
\end{Tm}
{\bf Proof:}\ \ {\bf I}. \ \ If (i) is true then (as $1\in A_0\,
\subset {\cal A} $)\ \ we have that\ \ $\delta_* (1) \in {\cal A}
\subset A_0^\prime$.\ \ In addition since\ \ $\delta$ \ \ is an
endomorphism of\ \ $\cal A$\ \ Proposition \ref{endom} tells us
that all the operators\ \ $\delta_*^k (1), \ k= 0, 1, ...$\ \ are
projections.\\ Thus (i) implies (ii).\\ The implication (ii)
$\Rightarrow$ (i) follows as a byproduct from Part {\bf II} of the
theorem (one can take ${\cal A} =\, _\infty A $).\ \ So let us
prove this part.\\

\noindent {\bf II}. \ \
 Note that in the situation under consideration all
 the conditions of Theorem \ref{An} are satisfied.
This implies the  commutativity of\ \    $_\infty A  $\ \ along
with (i) and (ii). The statement (iv) of
 Theorem \ref{An} also implies that\ \  $\delta$\ \ is an endomorphism of\ \ $_\infty A
 $\ \
and \ \ \mbox{$\delta \, :\, _nA \,  \to \, _{n-1}A, \ \ n=1,2,
... $.}\\ Clearly \ \ $\delta_*\,  :\,  _nA \, \to \,  _{n+1}A, \
\ n=0,1, ... $ \
\
 which proves (iii) along with the property  \ \
$\delta_* \, :\,  _\infty A \,  \hookrightarrow \, _\infty A $.\\
Since $$ UU^* = \delta(1) \in A_0 \subset\, _\infty A  \subset\,
(_\infty A )^\prime $$ Lemma \ref{V} implies that\ \ $\delta_*$\ \
is a morphism (and by the foregoing notes it is an endomorphism)
of\ \ $_\infty A $.\\ Finally as\ \ $UU^*, \ U^*U \in (_\infty A
)^\prime $\ \ Lemma \ref{V} implies (iv) as well.
\mbox{}\qed\mbox{}\\ The next result is the natural generalization
of the theorem just proved and it gives the solution to the
problem we are considering in this section.
\begin{Tm}
\label{A.inf}\ \ {\bf I}. \ \ The following statements are
equivalent:\\

\noindent (i) \ \ There exists a commutative $C^*-$algebra\ \
${\cal A } \supset A_0$\ \ such that both the mappings \ \
$\delta$\ \ and \ \ $\delta_*$ are endomorphisms of \ \ ${\cal
A}$.\\

\noindent (ii)\ \ All the operators \ \ $\delta_*^k (1), \ k= 0,
1, ...$\ \ are projections and $$ \delta_*^k (1) \in A_0^\prime ,
\ \ \ \delta^k (A_0) \subset A_0^\prime ,\ \ \ \delta_* (1) \in
[\delta^k (A_0)]^\prime , \ \  k= 0, 1, ... $$

{\bf II}. \ \ If condition (ii) of part {\bf I} is satisfied
then\\

\noindent 1)
\begin{equation}
\label{e.inf} _\infty (A_\infty )  = \{ A_\infty , \delta_*(
A_\infty   ), . . . , \delta_*^k ( A_\infty   ), . . . \}
\end{equation}
 is the minimal commutative $C^*-$algebra
containing $A_0$ and such that
\begin{equation}
\label{e.delta}
 \delta \, : \,
_\infty (A_\infty ) \, \to\,  _\infty (A_\infty ) \ \ \ \ \ and \
\ \ \ \ \delta_*\, :\,  _\infty (A_\infty )\,  \to\, _\infty
(A_\infty )
\end{equation}
 are endomorphisms.\\
Moreover we have: \\

\noindent (i) For any\ \  $0 \le l \le k$ \ \ $$ \delta_*^k (
A_\infty  )\cdot \delta_*^l ( A_\infty  )
 \subset \delta_*^k (A_\infty  )
$$

\noindent (ii) For every $n$ the subalgebra\ \ $_n(A_\infty )\,
\subset \,  _\infty (A_\infty ) $ \ \ given by (\ref{e.An}) (with\ \ $A_0$\ \ substituted by\ \ $A_\infty $)
 is the set of operators of the form
$$ \alpha_0 + \delta_*(\alpha_1) + ... + \delta_*^n(\alpha_n),\ \
\ \ \alpha_i \in  A_\infty  . $$

\noindent (iii)\ \
$
\delta \, : \,  _n(A_\infty )\,  \to \, _{n-1}(A_\infty ), \ \
n=1,2, ...  \ \ \
      \delta_* \, :\, _n(A_\infty )\,  \to \, _{n+1}(A_\infty ), \ \ n=0,1, ...
$\\
\ \\

\noindent (iv)\ \
\begin{equation}
\label{e.iv} U\alpha = \delta (\alpha )U \ \  and \ \ \ U^* \alpha
=\delta_*(\alpha )U^* \ \ \ \ for\   any \ \ \alpha \in\,  _\infty
(A_\infty ) .
\end{equation}\\

\noindent  2)
\begin{equation}
\label{e.Ainfty2} _\infty (A_\infty ) =\ ( _\infty A )_\infty
\end{equation}
where\ \ $ ( _\infty A )_\infty = \{ _\infty A  , \delta (_\infty
A ), . . . , \delta^k (_\infty A   ), . . . \}$\\ Moreover we
have:
\\

\noindent (i) For any\ \  $0 \le l \le k$ \ \ $$ \delta^k (
_\infty A  )\cdot \delta^l ( _\infty A  )
 \subset \delta^k (_\infty A  )
$$

\noindent (ii) For every $n$ the subalgebra \ \ $ ( _\infty A
)_n\, \subset \, ( _\infty A )_\infty  $\ \ given by (\ref{e.An1})
(with \ \ $A_0$ \ \ substituted by \ \ $_\infty A $)
 is the set of operators of the form
$$ \alpha_0 + \delta(\alpha_1) + ... + \delta^n(\alpha_n),\ \ \ \
\alpha_i \in  _\infty A  . $$

\noindent (iii)\ \
$
\delta_*\,  :\,   ( _\infty A )_n ,\ \to \, ( _\infty A )_{n-1}, \
\ n=1,2, ...  \ \ \ \ \       \delta \,  :\  ( _\infty A )_n \,
\to \, ( _\infty A )_{n+1}, \ \ n=0,1, ...  $\\
\end{Tm}
{\bf Proof:}\ \ {\bf I}. \ \ If (i) is true then  Proposition
\ref{endom} tells us that all the operators\ \ $\delta_*^k (1), \
k= 0, 1, ...$\ \ are projections.  In addition we have
$$\delta_*^k (1) \in {\cal A}\subset A_0^\prime ,$$ $$\delta^k
(A_0)\subset {\cal A} \subset A_0^\prime $$ and
 $$ \delta_* (1)
\in {\cal A}\subset [\delta^k (A_0)]^\prime , \ \ k= 0, 1, ... $$
 Thus (i) implies (ii).\\
As in the previous theorem the implication (ii) $\Rightarrow$ (i)
follows as a byproduct from Part {\bf II} of the theorem (one can
take ${\cal A} =\, _\infty (A_\infty ) $).\ \ So let us prove this
part.\\

\noindent {\bf II}. \ \
 Since in the situation under consideration all the
assumptions of \mbox{Theorem \ref{A1}} are satisfied it follows
that \ \ $A_\infty  $\ \ is a commutative algebra and
\begin{equation}
\label{e.i0} \delta \, :\,  A_\infty \, \hookrightarrow \,
A_\infty
\end{equation}
is an endomorphism and
\begin{equation}
\label{e.i1} \delta_* (1) \in \, (A_\infty  )^\prime
\end{equation}
In view of (\ref{e.i0}) and  (\ref{e.i1})  we can apply Theorems
\ref{A0.inf} and \ref{An}       (with $A_0$ substituted by
$A_\infty $) and thus all the statements in   1) are proved.\\ Let
us verify 2).\\
 Observe that in the situation considered all the assumptions of Theorem \ref{A2}
are satisfied (in particular \ \ $\delta (1) \in [\delta_*^k
(A_0)]^\prime , \ \  k= 0, 1, ...$\ \ since \ \ $\delta (1) \in \,
_\infty (A_\infty )$  \ \  and \ \ $\delta_*^k (A_0) \,  \subset
\,  _\infty (A_\infty )$). Thus
\begin{equation}
\label{e.i2} \delta_*\, :\,  _\infty A  \,  \hookrightarrow \,
_\infty A
\end{equation}
is an endomorphism.\\ Since $\delta (1)\, \in \,  _\infty
(A_\infty ) $ \ \ and \ \ $_\infty A \,  \subset \,  _\infty
(A_\infty )$\ \  it follows that
\begin{equation}
\label{e.i3} \delta (1) \in \,  (_\infty A  )^\prime
\end{equation}
In view of (\ref{e.i2}) and (\ref{e.i3}) we can apply Theorem
\ref{A0.inf} (where we exchange \mbox{$\delta$\ \  for\ \
$\delta_*$)} to the algebra\ \  $_\infty A $ \ \ (instead of \ \
$A_0$)
 and it follows that both the mappings
\begin{equation}
\label{e.00} \delta\,  :\, ( _\infty A )_\infty   \, \to \,  (
_\infty A )_\infty \ \ \ \ {\rm and} \ \ \ \ \delta_*\,  :\, (
_\infty A )_\infty  \, \to \,  ( _\infty A )_\infty
\end{equation}
are endomorphisms and the properties (i), (ii) and (iii) of 2) are
true as well.\\ Moreover (\ref{e.00}) along with  the construction
of \ \ $ _\infty (A_\infty )$ \ \ and\ \ $ ( _\infty A )_\infty$\
\  imply (\ref{e.Ainfty2}). \mbox{}\qed\mbox{}

\section{The structure of the algebra $\{ A_0, U  \}$ }
\label{3} \setcounter{equation}{0}

Throughout this section we fix   a partial
isometry $U$ and a commutative algebra  $A_0 \subset L(H)$ containing the identity and
 satisfying the assumptions (ii) of Part I of Theorem \ref{A.inf}.\\

\noindent The purpose of the section is to  describe the structure
of the $C^*-$algebra  $\{ A_0, U  \}$  generated by $A_0$ and
$U$.\\

\noindent To simplify the notation we denote by $A$ the algebra
\begin{equation}
\label{e.A} A = \, _\infty (A_\infty ) = \,   ( _\infty A )_\infty
\end{equation}
where \ \ $ _\infty (A_\infty )$\ \ and \ \ $( _\infty A
)_\infty$\ \ are the algebras defined by \ \ (\ref{e.inf})\ \ and
\ \ (\ref{e.Ainfty2})\ \ respectively.\\

\noindent As it will be clarified  the algebra \ \ $A$\ \ plays
here the crucial role of the 'coefficient' algebra.
\begin{La}
\label{3a.4}
 The algebra $\{ A_0, U  \}$  is the uniform closure
of  finite sums of the form
\begin{equation}
\label{e3a.4} U^{n*}\beta_{-n} + . . . + \beta_0 + . . . +
\beta_nU^n
\end{equation}
where \ \ $\beta_k \in \,  A$\ \ (\ref{e.A})\ \
 and\ \  $\beta_{\pm k}=
U^kU^{k*}\beta_{\pm k} = \beta_{\pm k}U^kU^{k*}, \ \ k=0,...,n$.\\
This algebra is also the uniform closure of  finite sums of the
form
\begin{equation}
\label{e3a.4'} \alpha_{-n}U^{n*} + . . . + \alpha_0 + . . . +
U^n\alpha_n
\end{equation}
where \ \ $\alpha_k \in\, A$\ \
 and \ \ $\alpha_{\pm k}=
U^{k*}U^k\alpha_{\pm k} = \alpha_{\pm k}U^{k*}U^k, \ \
k=0,...,n$.\\
\end{La}
{\bf Proof:} Follows in a routine way from the definition of $\{
A_0, U  \}$ along with (\ref{e.delta}), (\ref{e.iv}) and the
observation that $$ U^kU^{k*}, \ U^{k*}U^k \in\,
 _\infty
(A_\infty ) = A\  \ \ \ \ k=1,2, ... $$ and the equalities $$
U^{k*}U^kU^{k*} =U^{k*}, \ \ \ \ U^kU^{k*}U^k = U^k ,\ \ \ \ \
k=1,2, . . . $$
 \mbox{}\qed\mbox{}
\begin{I}
\label{*}
 \em We shall
say that  the algebra  $\{ A_0, U  \}$ possesses the {\em
property} $(*)$ if for every element $b\in \{ A_0, U  \}$ having
the form (\ref{e3a.4}) the following inequality holds:
\begin{equation}
\label{e3a.5} \Vert b \Vert \ge \Vert \beta_0 \Vert
\end{equation}
Clearly the algebra  $\{ A_0, U  \}$ possesses the  property $(*)$
iff for every element $b\in \{ A_0, U  \}$ having the form
(\ref{e3a.4'}) the following inequality holds:
\begin{equation}
\label{e3a.5'} \Vert b \Vert \ge \Vert \alpha_0 \Vert
\end{equation}
Indeed. Due to (\ref{e.iv}) one can always pass from (\ref{e3a.4}) to
(\ref{e3a.4'}) (and back) leaving\ \  $\beta_0  = \alpha_0 $\ \
invariant under the passage.
\end{I}
\begin{Tm}
\label{3a.4} If $\{ A_0, U  \}$ possesses the  property $(*)$ then
for any element $b$ of the form (\ref{e3a.4}) we have
\begin{equation}
\label{e3a.5} \Vert b \Vert \ge \mathop {\max }\limits_{i =
\overline{0,n}} \Vert  \beta_{\pm i} \Vert
\end{equation}
and for any element $b$ of the form (\ref{e3a.4'}) we have
\begin{equation}
\label{e3a.5'} \Vert b \Vert \ge \mathop {\max }\limits_{i =
\overline{0,n}} \Vert  \alpha_{\pm i} \Vert
\end{equation}
\end{Tm}
{\bf Proof:} Let us prove (\ref{e3a.5}) first.\\
Fix some $k\in \overline{1,n}$. Since $\Vert U^k \Vert \le 1$ it follows that
\begin{equation}
\label{e.3a4-1}
\Vert b \Vert \ge \Vert U^k b\Vert
\end{equation}
We have
\begin{equation}
\label{e.3a4-2}
\begin{array}{*{20}l}
U^k b = U^k U^{n*} \beta_{-n}+ ... + U^k U^{k*} \beta_{-k}+ ... +\\
\ \\
U^k U^{s*} \beta_{-s} + ... +
U^k \beta_0 + ... + U^k \beta_n U^n
\end{array}
\end{equation}
Note that\\

\noindent 1) For $k> l$ we have (using the equality\ \
$U^{s*}U^sU^{s*} = U^s$,\ \ (\ref{e.iv}) and bearing in mind that
$U^sU^{s*} \in\, _\infty (A_\infty ) = A$):
\begin{eqnarray}
\label{e.3a4-3} && U^kU^{l*}\beta_{-l} =
U^kU^{k*}U^{(l-k)*}\beta_{-l} = \delta^k (1)U^{(l-k)*}\beta_{-l} =
\nonumber \\ &&\ \nonumber \\ &&U^{(l-k)*}\delta^l (1)\beta_{-l} =
U^{(l-k)*}(U^{l-k}U^{(l-k)*}\delta^l (1)\beta_{-l})=
U^{(l-k)*}\beta^\prime_{-l+k}
\end{eqnarray}
where\ \  $\beta^\prime_{-l+k} =
U^{l-k}U^{(l-k)*}\beta^\prime_{-l+k}$\ \ and \ \
$\beta^\prime_{-l+k}  \in A$.\\

\noindent 2)
\begin{equation}
\label{e.3a4-4} U^k U^{k*} \beta_{-k} = \beta^\prime_0 \in A
\end{equation}

\noindent 3) For $k\ge l>0$ we have
\begin{eqnarray}
\label{e.3a4-5} &&U^k U^{l*}\beta_{-l} =
U^{k-l}U^lU^{l*}\beta_{-l} = \delta^{k-l} (U^lU^{l*}\beta_{-l}
)U^{k-l} = \nonumber \\ &&\ \nonumber
 \delta^{k-l} (U^lU^{l*}\beta_{-l} )U^{k-l}U^{(k-l)*}U^{k-l} =
\beta^\prime_{k-l}U^{k-l}
\end{eqnarray}
where\ \  $\beta^\prime_{k-l} =
\beta^\prime_{k-l}U^{k-l}U^{(k-l)*}$\ \ and \ \
$\beta^\prime_{k-l}  \in A$.\\

\noindent 4)
\begin{equation}
\label{e.3a4-6} U^k  \beta_0 = \delta^{k} (\beta_0  )U^k U^{k*}U^k
= \beta^\prime_k U^k
\end{equation}
where\ \  $\beta^\prime_{k} = \beta^\prime_{k}U^{k}U^{k*}$\ \ and
\ \ $\beta^\prime_{k}  \in A$.\\

\noindent 5) For $ l>0$ we have
\begin{equation}
\label{e.3a4-7} U^k \beta_lU^l = \delta^{k}(\beta_l ) U^{l+k}  =
\delta^{k}(\beta_l ) U^{l+k} U^{(l+k)*} U^{l+k} =
\beta^\prime_{l+k}U^{l+k}
\end{equation}
where\ \  $\beta^\prime_{l+k} =
\beta^\prime_{l+k}U^{l+k}U^{(l+k)*}$\ \ and \ \
$\beta^\prime_{l+k}  \in A$.\\

\noindent Now (\ref{e.3a4-1}), (\ref{e.3a4-2}) and (\ref{e.3a4-3})-(\ref{e.3a4-7})  along with the property $(*)$
imply
\begin{eqnarray}
\label{e.3a4-8}
&&\Vert b \Vert \ge \Vert U^k b \Vert = \Vert U^{(n-k)*}\beta^\prime_{n-k} + ... +  \beta^\prime_0
+ ... + \beta^\prime_{n+k}U^{n+k}\Vert \ge \nonumber \\
&&\ \nonumber \\
&&\Vert  \beta^\prime_0 \Vert = \Vert U^kU^{k*}\beta_{-k} \Vert = \Vert \beta_{-k}\Vert
\end{eqnarray}
Since\ \  $\Vert b \Vert = \Vert b^* \Vert$\ \ inequality  (\ref{e.3a4-8}) being applied to $b^*$ implies
$$
\Vert b \Vert=\Vert b^* \Vert \ge \Vert \beta^*_k \Vert = \Vert \beta_k \Vert
$$
thus finishing the proof of (\ref{e3a.5}).\\
(\ref{e3a.5'}) follows from (\ref{e3a.5}) by exchanging $U $ for $U^*$.
 \mbox{}\qed\mbox{}
\begin{I}
\label{N}
\em
 Thus if $\{ A_0, U  \}$ possesses
the property $(*)$ then all the coefficients in (\ref{e3a.4}) and
(\ref{e3a.4'}) are uniquely determined and\ \  $\alpha_0 = \beta_0$\ \
and for every $n=0, \pm 1, ...$ the mapping
\begin{equation}
\label{e3a.6} N_n : \{ A_0, U  \} \to A
\end{equation}
given by
\begin{equation}
\label{e3a.7} N_n (b)= \beta_n
\end{equation}
is correctly defined (as well as the corresponding mapping
exploiting the coefficients $\alpha_n$ is correctly defined).
\end{I}
The next Lemma \ref{3a.sum} gives a number of  norm estimates of sums of elements
in $C^*-$algebras. This lemma being useful in its own right also plays an important role in the proof of
 Theorem \ref{3a.N} that presents a certain formula for the norm calculation
of the elements of $ \{ A_0, U  \}$.\\
 The estimates presented in the lemma are probably
 known (in particular the components of the statement of the  lemma are given in
\cite{Ant}, Lemma 7.3 and \cite{AntLeb}, Lemma 22.3)). The proof
of the lemma can be obtained as a simple modification of the
reasoning given in the proof of \cite{AntLeb}, Lemma 22.3.\\
\begin{La}
\label{3a.sum}
For any $C^*-$ algebra $B$ and any elements $d_1, ..., d_m \in B$ we have
\begin{equation}
\label{be3.81}
\left\Vert  \sum_{i=1}^m d_i  \right\Vert^2 \le m \left\Vert  \sum_{i=1}^m d_i d_i^* \right\Vert
\end{equation}
and
\begin{equation}
\label{be3.81a}
\left\Vert  \sum_{i=1}^m d_i  \right\Vert^2 \le m \left\Vert  \sum_{i=1}^m d_i^* d_i \right\Vert
\end{equation}
On the other hand
\begin{equation}
\label{be3.82}
\left\Vert  \sum_{i=1}^m |d_i|  \right\Vert^2 \ge\frac{1}{m} \left\Vert  \sum_{i=1}^m d_i^* d_i \right\Vert
\end{equation}
and
\begin{equation}
\label{be3.82a}
\left\Vert  \sum_{i=1}^m \sqrt{d_i d_i^*}  \right\Vert^2 \ge\frac{1}{m} \left\Vert  \sum_{i=1}^m d_i d_i^* \right\Vert
\end{equation}
\end{La}
\begin{Tm}
\label{3a.N}
Let\ \   $\{ A_{0}, U  \}$\ \    be an algebra  such that  the pair \ \
$A_{0}, U $\ \  satisfies the assumptions (ii)of Part I of Theorem \ref{A.inf}.
If
the algebra\ \  $\{ A_{0}, U  \}$\ \  possesses
the property $(*)$ then for any element $b$ of the form (\ref{e3a.4})
we have
\begin{equation}
\label{be3.131}
\Vert b \Vert = \lim_{k\to\infty} \sqrt[4k]{ \left\Vert N_0 \left[ (bb^*)^{2k}\right]\right\Vert }
\end{equation}
where $N_0$ is the mapping defined by (\ref{e3a.7}).
\end{Tm}
{\bf Proof:} Applying (\ref{be3.81}) to the operator $$ b=
U^{n*}\beta_{-n} + . . . + \beta_0 + . . . + \beta_nU^n = d_{-n} +
... + d_0 + ... + d_n $$ we obtain $$ \Vert b \Vert^2 \le (2n+1)
\left\Vert  \sum_{i=-n}^n d_i d_i^* \right\Vert = (2n+1) \Vert N_0
(bb^*)\Vert $$ where $$ d_id_i^* = U^{-i*}\beta_i \beta_i^* U^{-i}
\ ,\ \ \  i<0; $$ $$ d_id_i^* = \beta_iU^iU^{i*}\beta_i^*\ , \ \ \
\ i\ge 0, $$ in either case $d_id_i^* \in\,  _\infty (A_\infty ) =
A $.\\ On the other hand as \ \  $\{ A_{0}, U  \}$\ \ possesses
the property $(*)$ we have $$ \Vert b \Vert^2 = \Vert bb^* \Vert
\ge \Vert N_0 (bb^*)\Vert $$ thus
\begin{equation}
\label{be3.101}
\Vert N_0 (bb^*)\Vert \le  \Vert bb^* \Vert = \Vert b \Vert^2  \le (2n+1) \Vert N_0 (bb^*)\Vert
\end{equation}
Applying (\ref{be3.101}) to $(bb^* )^k$ and having in mind that
$(bb^* )^k = (bb^* )^{k*} $ and $\Vert (bb^* )^{2k}\Vert = \Vert b \Vert^{4k} $ one has
\begin{equation}
\label{be3.111}
\Vert N_0 \left[ (bb^*)^{2k}\right]\Vert \le
 \Vert (bb^*)^k \cdot  (bb^*)^{k*}  \Vert =
\Vert b \Vert^{4k}  \le (4kn+1) \Vert N_0 \left[ (bb^*)^{2k}\right]\Vert
\end{equation}
since being written in the form (\ref{e3a.4}) $(bb^*)^k $ has not more than
$(4kn+1)$ summands.\\
So
\begin{equation}
\label{be3.121}
\sqrt[4k]{ \left\Vert N_0 \left[ (bb^*)^{2k}\right]\right\Vert }\le \Vert b \Vert \le
\sqrt[4k]{ 4kn+1 }\cdot \sqrt[4k]{  \left\Vert N_0 \left[ (bb^*)^{2k}\right]\right\Vert }
\end{equation}
Observing the equality
$$
\lim_{k\to\infty}\sqrt[4k]{ 4kn+1 } =1
$$
we conclude that
$$
\Vert b \Vert = \lim_{k\to\infty} \sqrt[4k]{ \left\Vert N_0 \left[ (bb^*)^{2k}\right]\right\Vert }
$$
The proof is complete.
\mbox{}\qed\mbox{}\\
The next result shows that the property $(*)$ plays the crucial
role in the determination of the structure of the algebra
 $\{ A_0, U  \}$
once the structure of\ \  $A =\, _\infty (A_\infty )$ \ \ is
elucidated.
\begin{Tm}
\label{3a.5} Let\ \   $\{ A_0^1, U_1  \}$\ \    and \ \ $\{ A_0^2,
U_2 \} $\  \ be two algebras such that both the pairs \ \ $A_0^i,
U_i ,\ \ i=1,2 $\ \  satisfy the assumptions (ii) of Part I of Theorem
\ref{A.inf}. Suppose that there exists an isomorphism $$ \varphi :
A_0^1 \to A_0^2 $$ such that under the mapping\ \  $U_1\to U_2$ \
\ the isomorphism\ \  $\varphi$ \ \ give rise to the isomorphism
\begin{equation}
\label{e3.''}
 \varphi :\, A^1
\to\, A^2
 \end{equation}
 (where\ \
$A^i =\, _\infty (A_\infty^i) , \ \ i=1,2$\ \
 are the  algebras given by (\ref{e.inf})
defined  by the pairs\ \   $A_0^i, U_i$\ \   respectively).\\ If
both the algebras\ \  $\{ A_0^i, U_i  \}, \ \ i=1,2$\ \  possess
the property $(*)$ then the mappings $$
 \varphi :  A_0^1 \to  A_0^2, \ \ \ U_1\to U_2
$$
give rise to the isomorphism
 $$
  \{ A_0^1, U_1 \}\cong \{ A_0^2, U_2  \}
  $$
\end{Tm}
{\bf Proof:} Consider an operator\ \ $b \in  \{ A_0^1, U_1 \}$\ \
having the form
\begin{equation}
\label{e3.01}
b  = U_1^{n*}\beta_{-n} + . . . +\beta_0 + . . . +
\beta_n U_1^n
\end{equation}
where $\beta_k \in\,  A^1 $.\\ Let
 \ \ $\varphi (b) \in \{ A_0^2, U_2 \}$ be the operator given by
\begin{equation}
\label{e3.02}
\varphi (b)  = U_2^{n*}\varphi (\beta_{-n}) + . . . +\varphi ( \beta_0) + . . . +
\varphi (\beta_n) U_2^n
\end{equation}
To prove the theorem it is enough to verify the equality\ \  $\Vert b \Vert = \Vert \varphi (b) \Vert$.\\
By Theorem \ref{3a.N} we have
\begin{equation}
\label{be3.013}
\Vert b \Vert = \lim_{k\to\infty} \sqrt[4k]{ \left\Vert N_0 \left[ (bb^*)^{2k}\right]\right\Vert }
\end{equation}
where
$$
 N_0 :  \{ A_0^1, U_1 \} \to\,  A^1
$$
is described in \ref{N}.\\
Similarly
\begin{equation}
\label{be3.013'}
\Vert \varphi (b) \Vert = \lim_{k\to\infty} \sqrt[4k]{ \left\Vert N_0 \left[ (\varphi (b) \varphi (b^*))^{2k}\right]\right\Vert }
\end{equation}
where
$$
 N_0 :  \{ A_0^2, U_2 \} \to A^2
$$
Observe that (\ref{e3.01}) and (\ref{e3.02}) along with the assumptions of the theorem imply
$$
 N_0 \left[ (\varphi (b) \varphi (b^*))^{2k}\right] =  N_0 \left[ (\varphi (bb^*))^{2k}\right] =
\varphi \left(  N_0 \left[ (bb^*)^{2k}\right] \right)
$$
and therefore (in view of (\ref{e3.''}))
$$
 \left\Vert N_0 \left[ (\varphi (b) \varphi (b^*))^{2k}\right] \right\Vert  =
 \left\Vert  N_0 \left[ (bb^*)^{2k}\right] \right\Vert
$$
This along with (\ref{be3.013'}) and  (\ref{be3.013}) implies the equality
$$
\Vert b \Vert = \Vert \varphi (b) \Vert
$$
and finishes the proof.
\mbox{}\qed\mbox{}\\
\begin{Rk}
\label{cross} \em It is worth mentioning that  the property $(*)$
and the results of Theorem \ref{3a.5} type play a fundamental role
in the theory of crossed products of $C^*-$algebras by discrete
groups (semigroups) of automorphisms (endomorphisms). Namely this
property is a {\em characteristic} property of the crossed product
and it enables one to construct its  faithful representations. The
importance of the property $(*)$ for the first time (probably) was
clarified  by O'Donovan \cite{O'Donov} in connection with the
description of  $C^*-$algebras generated by weighted shifts. The
most general result establishing the crucial role of this property
in the theory of crossed products of $C^*-$algebras by discrete
groups of {\em automorphisms} was obtained in \cite{Leb1} (see
also \cite{AntLeb}, Chapters 2,3 for complete proofs and various
applications) for an {\em arbitrary} $C^*-$algebra and {\em
amenable} discrete group. The relation of the corresponding
property to the faithful representations of crossed products by
{\em endomorphisms} generated by {\em isometries} was investigated
in \cite{BKR,ALNR}. The properties of this sort proved to be of
great value not only in pure $C^*-$theory but also in various
applications such as, for example, the construction of  symbolic
calculus and developing the solvability theory of functional
differential equations (see \cite{AnLebBel1, AnLebBel2}).
\end{Rk}

\section{The structure of the algebra ${\bf \cal B} = \{ 1, |a|, U \}$.}
\label{4} \setcounter{equation}{0}

Now we return to the consideration of  the initial object of our
investigation. Throughout this section we fix an operator $a\in
L(H)$ satisfying relation (\ref{eI.1}), take the partial isometry
$U$ from the polar decomposition (\ref{eI.2}) and set
 $A_0 = \{ 1,
|a| \} $.\\

\noindent  The aim of the section is to describe the structure of
the $C^*-$algebra ${\bf \cal B} = \{ 1, |a|, U \}$.\\

\noindent Following the path already exploited in the preceding
sections we start with the examination of the properties of the
partial isometry $U$ along with the properties of the induced
mappings and the arising algebras (extensions of $A_0$).\\

\begin{Tm}
\label{2.2}
Let $U$ be the partial isometry defined by  polar decomposition (\ref{eI.2})
of an operator $a$ satisfying  relation (\ref{eI.1}) then\\

\noindent 1). $U^*U \in \{ 1, |a| \}^{\prime \prime}$,\\ where we
denote by\ \ ${\cal A}^{\prime \prime}$\ \ the bicommutant of an
algebra\ \ $\cal A$\ \ (that is the Von Neumann algebra generated
by\ \ $\cal A$).\\

\noindent 2). $U^k U^{*k} \in \{1, |a| \}^{\prime \prime}, \ \ k=1,2, ... $\\

\noindent 3). $[U^{*l}U^l , U^k U^{*k}] =0 \ , \ \ k=1,2, ...; \ \ l=1,2, ...  $\\

\noindent 4). $U^* U^k U^{*l} = U^{k-1}U^{*l}, \ \ 1\le k\le l $\\

\noindent 5). $(U^k U^{*k})^2 = U^k U^{*k},$ \ that is $U^k$ and $U^{*k}$ are partial isometries.\\

\noindent 6). If $k\ge l$ then
$$
U^{*k} U^k U^{*l} U^l =   U^{*l} U^l U^{*k} U^k = U^{*k} U^k
$$
and
$$
U^k U^{*k}U^l U^{*l} = U^l U^{*l}U^k U^{*k} = U^k U^{*k}
$$\\

\noindent 7).
\begin{equation}
\label{01} \delta : \{ |a| \} \to  \delta(\{  |a| \})\subset  \{
1, |a| \}
\end{equation}
is a morphism and
\begin{equation}
\label{I.b}
 U b = \delta (b)U \ \ {\rm  and} \ \   U^*  \delta (b) =bU^*, \ \  b \in \{ |a| \}
\end{equation}




\noindent 8). ${\delta}^k \{ |a|\} \subset \{ 1, |a|, UU^*, ... ,
U^{k-1}U^{k-1*}  \}, \ \ \ k=1,2, ... $\\

\noindent 9). ${\delta}^k (\{ |a|\} ) = U^k U^{*k}{\delta}^k (\{
|a| \} ) = {\delta}^k  (\{ |a|\} )  U^k U^{*k}$\\ that is $U^k
U^{*k}$ is the projection onto the invariant subspace where the
whole of the algebra ${\delta}^k  (\{ |a| \} ) $ acts.\\

\noindent 10). For every\ \  $\alpha \in \{ |a| \}$\ \ we have\ \
$\delta_*\delta(\alpha ) = \alpha = U^*U\alpha  = \alpha \, U^*U$.
\end{Tm}
{\bf Proof:}\\
1).
 We have that $U^*U$ is the orthogonal projection onto
${\rm Im}\, |a| = ({\rm Ker}\, |a|)^\perp$. Therefore $U^*U$ is
the spectral projection of $|a|$ corresponding to the set\ \
$\sigma (|a|) \setminus \{ 0 \}$ \ \ (where we denote by $\sigma
(\alpha )$ the spectrum of an operator $\alpha$). Thus by the
spectral theorem
 (see, for example, \cite{Naimark}, \S 17)\ \  $U^*U \in \{ 1, |a| \}^{\prime \prime}$.\\

\noindent 7). The statement that \ \ $ \delta : \{ |a| \} \to
\delta(\{  |a| \})$\ \ is a morphism and (\ref{I.b})  follow from
1) in view of part I of  Lemma \ref{V}. In addition (\ref{I.a})
means that \ \ $\delta (|a|)\in  \{ 1, |a| \}$\ \  and therefore
 $ \delta  \{ |a| \} \subset  \{ 1, |a| \}$.\\


\noindent 8). By (\ref{01}) we have that $$ {\delta} \{ |a|\}
\subset \{ 1, |a| \} $$ Thus $$ {\delta}^2 \{ |a|\} \subset
{\delta}(\{ 1, |a| \}) = \{ {\delta}(1), {\delta}(|a|)\} \subset
\{ UU^*, 1, |a|\}. $$ Continuing this reasoning we obtain 8).\\

\noindent 2). Since $U^*U$ is the spectral projection
corresponding to the interval $(0, |a|]$ (see the proof of 1)) it
follows that there exists a sequence $\alpha_n$ of elements of $\{
|a|  \}$ such that
\begin{equation}
\label{strong}
\alpha_n\ \  \mathop{\longrightarrow}\limits^{strongly}\ \ U^*U
\end{equation}
Therefore we have $$ U\alpha_n U^* \ \
\mathop{\longrightarrow}\limits^{strongly}\ \ U (U^*U)U^* = UU^*
$$ Since (by 8)) $U\alpha_n U^* \subset \{ 1,|a| \}$ it follows
that $UU^* \subset  \{ 1,|a| \}^{\prime\prime}$.\\ The further
proof goes by induction.\\ Suppose that\ \  $U^kU^{k*}\in  \{
1,|a| \}^{\prime\prime} , \ \ k=\overline{1, n-1}$. Taking the
sequence $\alpha_n$ mentioned above we have
\begin{equation}
\label{e2.2.1}
U^n\alpha_n U^{n*} \ \  \mathop{\longrightarrow}\limits^{strongly}\ \ U^{n} (U^*U)U^{n*} =
U^{n-1} (UU^*U) U^{n*} =
 U^{n}U^{n*}
\end{equation}
But due to 8) and the assumption of the induction we have $$
U^n\alpha_n U^{n*} \subset \{ 1, |a|, UU^*, ... , U^{n-1}U^{n-1*}
\}\subset  \{ 1,|a| \}^{\prime\prime} $$ and therefore
(\ref{e2.2.1}) implies \ \ $U^{n}U^{n*} \in  \{ 1,|a|
\}^{\prime\prime}$. So 2) is proved.\\

\noindent 3), \ 4), \ 5) \  and 6).\  These follow from 1) and 2) along with Lemma \ref{""}.\\

\noindent 9). In view of 5) we have  the relations \ \
$U^kU^{*k}U^k = U^k$ \ \  and \ \ $U^{*k}U^kU^{*k} =U^{*k}$\ \
that imply 9).\\

\noindent 10). Recall that since $U^*U$ is the orthogonal
projection onto ${\rm Im}\, |a|=  ({\rm Ker}\, |a|)^\perp$ it
follows that\ \ $U^*U |a| =|a|= |a| U^*U$. This along with 1)
implies 10). \mbox{}\qed\mbox{}
\begin{Rk}
\label{.} \em Most of the properties listed in Theorem \ref{2.2}
are known (we have presented them for the sake of completeness).
In particular, one can find some generalizations of the properties
3) and 7)  in Propositions 28 and 29 in \cite{OstSam}, Section 2.1
that also contains a lot of important information on the subject.
\end{Rk}
\begin{I}
\label{}
\em \
Now  as in section \ref{1}
we introduce the necessary for our future goals
  extensions  of the algebra $A_0=\{1, |a|  \}$.\\
Having in mind \ref{defAn} and statements 2) and 6) of the theorem just proved
for every  $n= 0,1, ...$ we set
\begin{equation}
\label{eA-n}
 A_n  =\{ 1, |a|, UU^*, ... , U^nU^{n*}  \} \subset \{ 1, |a| \}^{\prime \prime}
 \end{equation}
Clearly $$ \delta ( A_n)\subset  {A}_{n+1} $$ is a morphism (the
latter follows from Lemma \ref{V}).\\ Now bearing in mind
(\ref{e.A1}) we set
\begin{equation}
\label{eA-infty} A_\infty  = \{ 1, |a|, UU^*, ... , U^nU^{n*}, ...
\}\subset \{ 1, |a| \}^{\prime \prime}
\end{equation}
Evidently
\begin{equation}
\label{e.s} \delta : A_\infty  \to A_\infty  \ \ \mbox{\rm is an
endomorphism}.
\end{equation}
and in particular
\begin{equation}
\label{e.s'}
 \delta^k (A_0)\subset  A_\infty  \subset
 {A_0}^\prime
\end{equation}
 Note also that by
statement 1) of Theorem \ref{2.2} and (\ref{eA-infty}) we have
\begin{equation}
\label{e.s1} \delta_* (1)=U^{*}U \in  (A_\infty  )^{\prime}
\subset [ \delta^k (A_0)]^\prime , \ \ \ k=0,1,...
\end{equation}
which along with (\ref{e.s}) implies (by Lemma \ref{endom1})
\begin{equation}
\label{UkU*k} \delta_*^k (1)=U^{*k}U^k \in  (A_\infty )^{\prime}
\subset A_0^\prime \ , \ \ k=1,2, ...
\end{equation}
Moreover (\ref{e.s}) along with (\ref{e.s1}) also imply (by Lemma \ref{V})
\begin{equation}
\label{e3.1}
 U b = \delta(b)U \ \  {\rm and}\ \
U^*\delta(b)= bU^*, \ \ \ b\in A_\infty
\end{equation}\\

\noindent Observe now that (\ref{e.s'}), (\ref{e.s1}) and
(\ref{UkU*k})
 mean
that for the algebra $A_0 = \{ 1, |a| \}$ and
 the operator $U$ under consideration
 all the assumptions of Theorem  \ref{A.inf} are satisfied.\\

\noindent  Thus we arrived at the next statement.
\end{I}
\begin{Tm}
\label{3.3'} Let $a$ be an operator satisfying relation
(\ref{eI.1}) and  $U$ be the partial isometry defined by  polar
decomposition (\ref{eI.2}). Let $A_\infty $ be the algebra given
by (\ref{eA-infty})
 then\\

 \noindent 1)
\begin{equation}
\label{e.inf'}
 A =\, _\infty (A_\infty )  = \{
A_\infty   , \delta_*( A_\infty   ), . . . , \delta_*^k ( A_\infty
), . . .  \}
\end{equation}
 is the minimal commutative $C^*-$algebra
containing $A_0= \{ 1, |a|\}$ and such that
\begin{equation}
\label{e.delta'}
 \delta :
A \to A \ \ \ \ \ and \ \ \ \ \ \delta_* : A \to A
\end{equation}
 are endomorphisms.\\
Moreover we have: \\

\noindent (i) For any\ \  $0 \le l \le k$ \ \ $$ \delta_*^k (
A_\infty )\cdot \delta_*^l ( A_\infty  )
 \subset \delta_*^k (A_\infty  )
$$

\noindent (ii) For every $n$ the subalgebra
\begin{equation}
\label{e.A-infty-n}
 _n(A_\infty )  = \{ A_\infty  , ...,
\delta_*^n (A_\infty  ) \} \subset A
\end{equation}
  is the set of operators of the form
$$ \alpha_0 + \delta_*(\alpha_1) + ... + \delta_*^n(\alpha_n),\ \
\ \ \alpha_i \in  A_\infty  . $$

\noindent (iii)\ \
$
\delta :\,   _n(A_\infty )  \to\,  _{n-1}(A_\infty ), \ \ n=1,2,
... \ \ \ \ \
      \delta_* : \, _n(A_\infty )  \to \, _{n+1}(A_\infty ), \ \ n=0,1, ...
$\\ \ \\

\noindent (iv)\ \
\begin{equation}
\label{e.iv'} U\alpha = \delta (\alpha )U \ \  and \ \ \ U^*
\alpha =\delta_*(\alpha )U^* \ \ \ \ for\   any \ \ \alpha \in A.
\end{equation}

 \noindent  2) Let
 \begin{equation}
 \label{e.inftyA}
 _\infty A  = \{ A_0, \delta_*(A_0), ...,
 \delta_*^n(A_0),... \}
 \end{equation}
Then
\begin{equation}
\label{e.Ainfty2'} A=\, _\infty (A_\infty ) =\, ( _\infty A
)_\infty
\end{equation}
where $ ( _\infty A )_\infty = \{ _\infty A  , \delta (_\infty A
), . . . , \delta^k (_\infty A   ), . . . \}$\\ Moreover we have:
\\

\noindent (i) For any\ \  $0 \le l \le k$ \ \ $$ \delta^k (
_\infty A  )\cdot \delta^l ( _\infty A  )
 \subset \delta^k (_\infty A  )
$$

\noindent (ii)  For every $n$ the subalgebra
\begin{equation}
\label{e.A-inf--n} ( _\infty A )_n =
 \{ _\infty A  , ..., \delta^n (_\infty A  ) \}
 \  \subset \, A
\end{equation}
 is the set of operators of the form
$$ \alpha_0 + \delta(\alpha_1) + ... + \delta^n(\alpha_n),\ \ \ \
\alpha_i \in\,   _\infty A  . $$

\noindent (iii)\ \
$
\delta_*\  :\,    ( _\infty A)_n \, \to \, ( _\infty A )_{n-1}, \
\ n=1,2, ...  \ \ \ \ \
      \delta\  :\  ( _\infty A )_n\,   \to\,
 ( _\infty A )_{n+1}, \ \ n=0,1, ...  $\\
\end{Tm}
Now we proceed to the description of the algebra  ${\bf \cal B} =
\{ 1, |a|, U \}$.\\

\noindent In this situation Lemma \ref{3a.4} turns into
\begin{La}
\label{3.4}
The algebra ${\bf \cal B} = \{ 1, |a|, U \}$ is
 the uniform closure of  finite sums of the form
\begin{equation}
\label{e3.4}
U^{n*}\beta_{-n} + . . . + \beta_0 + . . . + \beta_nU^n
\end{equation}
where\ \ $\beta_k \in A$\ \  ($A$\ \ is given by (\ref{e.inf'}))\
\ and\ \ $\beta_{\pm k}= U^kU^{k*}\beta_{\pm k} = \beta_{\pm
k}U^kU^{k*}, \ \ k=0,...,n$.\\ This algebra is also the uniform
closure of finite sums of the form
\begin{equation}
\label{e3.4'} \alpha_{-n}U^{n*} + . . . + \alpha_0 + . . . +
U^n\alpha_n
\end{equation}
where \ \ $\alpha_k \in A$\ \
 and \ \ $\alpha_{\pm k}=
U^{k*}U^k\alpha_{\pm k} = \alpha_{\pm k}U^{k*}U^k, \ \
k=0,...,n$.\\
\end{La}
\begin{I}
\label{*'} \em
 Using the notation of   \ref{*} we say that the
algebra ${\bf \cal B}$ possesses the {\em property} $(*)$ if for
every element $b\in {\bf \cal B}$ having the form (\ref{e3.4}) the
following inequality holds:
\begin{equation}
\label{e3.5}
\Vert b \Vert \ge \Vert \beta_0 \Vert
\end{equation}
In view of Theorem \ref{3a.4} we conclude that \\

{ \em if ${\bf \cal B}$ possesses the  property $(*)$ then for any
element $b$ of the form (\ref{e3.4}) we have
\begin{equation}
\label{e3.5}
\Vert b \Vert \ge \mathop {\max }\limits_{i = \overline{0,n}}
\Vert  \beta_{\pm i} \Vert
\end{equation}
and for any element $b$ of the form (\ref{e3.4'}) we have
\begin{equation}
\label{e3.5'} \Vert b \Vert \ge \mathop {\max }\limits_{i =
\overline{0,n}} \Vert  \alpha_{\pm i} \Vert
\end{equation}
}
\end{I}

\noindent  In the situation considered  Theorem \ref{3a.5} turns
into
\begin{Tm}
\label{3.5}
Let $a_1 =U_1 |a_1|$ and $a_2 =U_2 |a_2|$ be the polar decompositions of  operators
$a_1$ and $a_2$ and both the operators $a_1$ and $a_2$ satisfy the condition (\ref{eI.1}).
Suppose that  the mapping
\begin{equation}
\label{e3}
|a_1| \to |a_2| , \ \ U_1 \to U_2
\end{equation}
gives rise to the isomorphism
 $$ A^1 \cong
 A^2 $$
 (where $A^i, \ \ i=1,2$
are defined by (\ref{e.inf'}) for the operators $a_i$ and $U_i$
respectively).\\ If both the algebras ${\bf \cal B}_i = \{1,
|a_i|, U_i  \}, \ \ i=1,2$ possess the property $(*)$ then the
mapping (\ref{e3}) gives rise to the isomorphism $$ {\bf \cal
B}_1\cong {\bf \cal B}_2 $$
\end{Tm}

\section{The structure of the algebra ${\bf B} = \{ 1, a \}$.}
\label{5} \setcounter{equation}{0}

Now we return to the consideration of the structure of the algebra
${\bf B} = \{ 1, a \}$. In this situation along with the mapping
$\delta$ (\ref{e.delta}) we shall introduce  the mapping $$
\varphi : \{ 1, |a|\} \to \{ 1,  \delta (|a|) \} \subset  \{ 1,
|a|\} $$ given by
\begin{equation}
\label{eb4.1} \varphi \left[ (a^*a)^k\right] =  (aa^*)^k = \delta
\left[ (aa^*)^k\right] =
 \delta^k (a^*a), \ \ \ \ \  \varphi (1)=1.
\end{equation}
One can easily verify that for every $b\in \{ 1, |a|\}$ we have
\begin{equation}
\label{eb4.2}
a b = \varphi (b) a \ \ \ {\rm and}\ \ \ ba^* = a^*
\varphi (b)
\end{equation}
Indeed. The first equality follows from the equalities
 $$
a\cdot (a^*a)^k = (aa^*)^k \cdot a = \varphi [(a^*a)^k]\cdot a $$
and $$ a\cdot 1 =1\cdot a =\varphi (1)\cdot a .$$ And the second
follows from the first one by passage to the adjoint operator.\\
 Note that for any $k,l \ge 0$
we also have the following relations:\\

if $k>l$ then
\begin{equation}
\label{eb4.3}
a^ka^{*l} = b a^{k-l} \ \  {\rm where}\ \ b\in  \{ 1, |a|\};
\end{equation}

if $k<l$ then
\begin{equation}
\label{eb4.4}
a^ka^{*l} =  a^{*(k-l)}b \ \  {\rm where}\ \ b\in  \{ 1, |a|\};
\end{equation}

and
\begin{equation}
\label{eb4.5}
a^ka^{*k} \in  \{ 1, |a|\}.
\end{equation}
where we take $a^0 = a^{*0} =1$.\\ Indeed. If  $k\ge l$ then using
(\ref{eb4.2}) several times we have
 $$
\begin{array}{*{20}l}
a^ka^{*l}= a^{k-1} (aa^*) a^{*l-1} = \varphi^{k-1}(aa^*)\cdot
a^{k-1}\cdot a^{*l-1} = \\ \ \\ \varphi^{k-1}(aa^*)\cdot
\varphi^{k-2}(aa^*)\cdot ... \cdot \varphi^{k-l}(aa^*)\cdot
a^{k-l} = b a^{k-l}
\end{array}
$$ where $ \varphi^{k-1}(aa^*)\cdot \varphi^{k-2}(aa^*)\cdot ...
\cdot \varphi^{k-l}(aa^*) = b \in \ \{ 1, |a| \}$. Thus
(\ref{eb4.3}) and (\ref{eb4.5}) are proved.\\

\noindent (\ref{eb4.4}) follows from (\ref{eb4.3}) since for $k<l$
we have $$
 a^ka^{*l} = \left(a^ka^{*l}\right)^{**} = \left(
a^la^{k*}\right)^* = \left(ba^{l-k}\right)^* = a^{*(l-k)}b^*  $$\\
 The next useful observation is given by
\begin{La}
\label{b1}
For any two operators of the form
$$
a^{*l_i}b_{l_i , m_i}a^{m_i}, \ \ \ i=1,2
$$
where $l_i , m_i \ge 0, \ \ b_{l_i , m_i} \in  \{ 1, |a|\} $ we have
$$
a^{*l_1}b_{l_1 , m_1}a^{m_1}\cdot a^{*l_2}b_{l_2 , m_2}a^{m_2} = a^{*l_3}b_{l_3 , m_3}a^{m_3}
$$
where $l_3 , m_3 \ge 0, \ \ b_{l_3 , m_3} \in  \{ 1, |a|\} $ and
$$
m_1 - l_1 + m_2 -l_2 = m_3 - l_3
$$
\end{La}
{\bf Proof:} Follows in a routine way from (\ref{eb4.2}) -
(\ref{eb4.5}). \mbox{}\qed\mbox{}\\
 Applying this lemma and  using again (\ref{eb4.2}) -
(\ref{eb4.5}) we conclude that  the algebra ${\bf B} = \{ 1, a \}$
is the uniform closure of finite sums of  elements of the form
\begin{equation}
\label{eb4.6} b= b_{-n} + ... + b_0 + ... + b_n , \ \ \  n\ge0
\end{equation}
where $b_k, \ \  -n\le k \le n$ is a finite sum of the form $$ b_k
= \sum_{m-l=k}a^{*l}b_{l,m}a^m , \ \ \ \  b_{l,m} \in \{ 1,
|a|\},\ \ l,m \ge 0 $$
\begin{I}
\label{B_0}
\em
Let ${\bf B}_0 $ be the subalgebra of ${\bf B} $ generated by the elements of the form
$$
a^{*l}b a^l, \ \ l\ge 0,\ \ b\in \{ 1, |a|\}
$$
Lemma \ref{b1}
 shows that  ${\bf B}_0 $ is indeed a $C^*-$subalgebra
of ${\bf B} $. Note also that ${\bf B}_0 $ is a commutative
algebra since it is a subalgebra of\ \ $_\infty (A_\infty )$\ \
(see (\ref{e.inf'})).\\
 Using the elements of ${\bf B}_0 $ one can rewrite
(\ref{eb4.6}) in the form
\begin{equation}
\label{0} b = a^{*n}\beta_{-n}+ ...+ \beta_0+ ... + \beta_n a^n
\end{equation}
where $\beta_k \in {\bf B}_0 , \ \ -n\le k\le n$ (cf
(\ref{e3.4})).
\end{I}
\begin{I}
\label{b1*} \em We shall say that the algebra ${\bf B} $ possesses
the {\em property $(**)$} if for any element $b\in {\bf B} $ of
the form (\ref{eb4.6}) the following inequality holds
\begin{equation}
\label{b1**}
\Vert b \Vert \ge \Vert b_0 \Vert
\end{equation}\\
{\em Remark.} One can easily see that if the algebra  ${\bf \cal B} = \{ 1, |a|, U \}$
 considered in the previous
section possesses the property $(*)$ (see (\ref{e3.5})) then the algebra ${\bf B} $
being a subalgebra of
${\bf \cal B}$ possesses the property $(**)$.\\
If  ${\bf B} $ possesses the property $(**)$ then the mapping
\begin{equation}
\label{bN_0}
{\bf N}_0 : {\bf B} \to {\bf B}_0
\end{equation}
given by
$$
{\bf N}_0 (b) =b_0
$$
for any $b$ having the form (\ref{eb4.6}) is correctly defined.
\end{I}
Now the analogue to Theorem \ref{3a.N} for the situation
considered is
\begin{Tm}
\label{3a.N'} Let $a$ be an operator satisfying relation
(\ref{eI.1})
  and\ \
${\bf B} = \{ 1, a \}$. If\ \ $\bf B$\ \ possesses the property
$(**)$ (\ref{b1**}) then then for any element $b$ of the form
(\ref{eb4.6}) we have
\begin{equation}
\label{be3.131'} \Vert b \Vert = \lim_{k\to\infty} \sqrt[4k]{
\left\Vert N_0 \left[ (bb^*)^{2k}\right]\right\Vert }
\end{equation}
where $N_0$ is the mapping defined by (\ref{bN_0}).
\end{Tm}
{\bf Proof:} The same as for Theorem \ref{3a.N}. One should apply
the estimates from Lemma \ref{3a.sum} to  operator (\ref{eb4.6})
$$
 b= b_{-n} + ... + b_0 + ... + b_n = d_{-n} + ... + d_0 + ...
 + d_n
 $$
 and observe that in view of Lemma \ref{b1} \ \
  $d_id^*_i \in {\bf B}_0 \ \ -n\le i\le n$.
  \mbox{}\qed\mbox{}\\

\noindent Following the same path one can obtain the next result
which is the analogue to Theorem \ref{3.5} for the algebra  under
consideration.
\begin{Tm}
\label{B1} Let\ \ ${\bf B}_i = \{ 1, a_i \}, \ \ i=1,2 $ where
both the operators $a_i , \ i=1,2$ satisfy the condition
(\ref{eI.1}). Suppose that the mapping $a_1 \mapsto a_2$ generates
the isomorphism $$ {\bf B}_{01} \cong {\bf B}_{02} $$ where ${\bf
B}_{0i}, \ i=1,2$ is the algebra described in \ref{B_0}
(corresponding to the algebra ${\bf B}_i$).\\ If both the algebras
${\bf B}_i , \ i=1,2$ possess the property $(**)$ then the mapping
$a_1 \mapsto a_2$ gives rise to the isomorphism $$ {\bf B}_{1}
\cong {\bf B}_{2} . $$
\end{Tm}

\noindent One can consider the results presented above from a bit
different point of view.
\begin{I}
\label{deg} \em Let us introduce a certain notion of 'degree' for
finite products of elements $a$ and $a^*$.\\ For any $l\ge 0$ we
set
\begin{equation}
\label{d}
{\rm deg}\, a^l =l \ \  {\rm and} \ \ {\rm deg}\, a^{*l} = -l .
\end{equation}
where we take $a^0 = a^{*0} =1$ and thus $ {\rm deg}\, 1=0$.\\
For a finite product $b$ of these  elements we set its degree ${\rm deg}\,b$
to be equal to the sum of the degrees of its factors.\\
  {\em Remark.} Of course in general   the number ${\rm deg}\,b$ depends {\em not} on the
element $b$ as an element of $L(H)$ but {\em on the representation} of $b$ in the form of a finite
product
of elements $a$ and $a^*$ and therefore it is reasonable to consider
the number $\deg$ as being defined on the free semigroup generated by $a$ and $a^*$.\\

\noindent Clearly ${\rm deg}\,b^* = - {\rm deg}\,b$ and  if ${\rm
deg}\,b_1 =k_1$, ${\rm deg}\,b_2 =k_2$
 then ${\rm deg}\,(b_1\cdot b_2) =k_1 +k_2$.\\

\noindent Evidently the algebra ${\bf B} = \{ 1, a\}$ is the
uniform closure of finite sums of the form
\begin{equation}
\label{b} b = b_{-n} + ... + b_0 + ... + b_n
\end{equation}
where $n \ge 0 $ and $b_i , \ \ -n\le i \le n$ is a finite sum of
the form $$ b_i = \sum \beta_{ki}, \ \ \ {\rm deg}\, \beta_{ki}=i
$$ One can verify that \\ {\em if an element $b$ is a finite
product of elements $a$ and $a^*$ and ${\rm deg}\,b =k$ then $b$
can be represented in the form $$ b = a^{*l}b_{l,m}a^m $$ where
$l,m \ge0, \ \ b_{l,m}\in \{ 1, |a| \}$ \ and \  $m-l=k$.}\\ So we
can reformulate      the notions from \ref{b1*} in the following
way.
\end{I}
\begin{I}
\label{b*} \em Let $\overline{\bf B}_0 $ be the subalgebra of
${\bf B} $ generated by finite products $b$ such that ${\rm deg}\,
b =0$.\\ Observe that $\overline{\bf B}_0 = {\bf B}_0$ where the
latter algebra is that defined in \ref{B_0}.\\

\noindent Thus the algebra ${\bf B} $ possesses the  property
$(**)$ iff for any element $b\in {\bf B} $ of the form (\ref{b})
the following inequality holds
\begin{equation}
\label{b**}
\Vert b \Vert \ge \Vert b_0 \Vert
\end{equation}
Therefore Theorem \ref{B1} is equivalent to the next statement.
\end{I}
\begin{Tm}
\label{B}
Let ${\bf B}_i = \{ 1, a_i \}, \ \ i=1,2 $  where both the operators
$a_i , \ i=1,2$ satisfy the condition (\ref{eI.1}). Suppose that the
mapping $a_1 \mapsto a_2$ generates
the isomorphism
$$
\overline{\bf B}_{01} \cong \overline{\bf B}_{02}
$$
where $\overline{\bf B}_{0i}, \ i=1,2$ is the algebra described in \ref{b*}
(corresponding to the algebra
${\bf B}_i$).\\
If for any element  $b_i \in {\bf B}_i $ having the form  (\ref{b}) the
inequality (\ref{b**}) holds  then the mapping
$a_1 \mapsto a_2$ gives rise to the isomorphism
$$
{\bf B}_{1} \cong {\bf B}_{2} .
$$
\end{Tm}


\begin{thebibliography}{                                  }


\bibitem{Kar}
V.P. Karassiov, {\em $sl(2)$ variational schemes for solving one
class of nonlinear quantum models,} Physics Letters A 238, 1998,
p. 19-28.


\bibitem{OdHorTer}

A. Odzijewicz, M. Horowski, A. Tereszkiewicz, {\em The integrable
multiboson systems and orthogonal polynomials,} IFT UwB (Preprint)
/14/ 2000.

\bibitem{MaxOdz}
V. Maximov, A. Odzijewicz, {\em The $q-$deformation of quantum
mechanics of one degree of freedom,} J. Math. Phys. {\bf 36},
1995, No 4, p. 1681-1690.

\bibitem{Odz1}
A. Odzijewicz, {\em Quantum Algebras and $q-$Special functions
Related to Coherent States Maps of the Disc,} Commun. Math. Phys.,
{\bf 192}, 1998, p. 183-215.



\bibitem{OstSam}
V. Ostrovskyi, Yu. Samoilenko, {\em Introduction to the Theory of
Representations of Finitely Presented $^*-$Algebras. I.
Representations by bounded operators,} Gordon and Breach, 1999.

\bibitem{Halm}
P.R. Halmos, {\em A Hilbert space problem book.} Van Nostrand, 1967.

\bibitem{Naimark}
M.A. Naimark, {\em Normed rings.} Nauka, Moscow, 1968 (Russian).

\bibitem{Ant}
A.B. Antonevich, {\em Linear functional equations. Operator
approach.} Universitetskoe Publishers, Minsk, 1988 (Russian).


\bibitem{AntLeb}
A. Antonevich, A. Lebedev, {\em Functional differential equations: I. $C^*-$theory.}
Longman Scientific $\&$ Technical, Pitman Monographs and Surveys in Pure and Applied
Mathematics 70, 1994.

\bibitem{O'Donov}
D.P. O'Donovan, {\em Weighted shifts and covariance algebras,}
Trans. Amer. Math. Soc. {\bf 208}, 1975, p. 1-25.

\bibitem{Leb1}
A.V. Lebedev, {\em On certain $C^*-$methods that are used while
investigating algebras associated with automorphisms and
endomorphisms,} Dep. VINITI, 1987, No 5351-B87 (Russian).

\bibitem{BKR}
S. Boyd, N. Keswani, I. Raeburn, {\em Faithful representations of
crossed products by endomorphisms,} Proceedings of the Amer. Math.
Soc. , {\bf 118}, 1993, No 2, p. 427-436.

\bibitem{ALNR}
S. Adji, M. Laca, M. Nilsen, I. Raeburn, {\em Crossed products by
semigroups of endomorphisms and the Toeplitz algebras of ordered
groups,} Proceedings of the Amer. Math. Soc. , {\bf 122}, 1994, No
4, p. 1133-1141.

\bibitem{AnLebBel1}
A. Antonevich, M. Belousov, A. Lebedev, {\em Functional
differential equations: II. $C^*-$applications. Part 1 Equations
with continuous coefficients,} Addison Wesley Longman, Pitman
Monographs and Surveys in Pure and Applied Mathematics 94, 1998.


\bibitem{AnLebBel2}
A. Antonevich, M. Belousov, A. Lebedev, {\em Functional
differential equations: II. $C^*-$applications. Part 2 Equations
with discontinuous coefficients and boundary value problems,}
Addison Wesley Longman, Pitman Monographs and Surveys in Pure and
Applied Mathematics 95, 1998.

\end{thebibliography}
\end{document}